\documentclass[12pt]{article}
\usepackage{amssymb}
\usepackage{amsthm}
\usepackage{graphicx}
\usepackage{graphics}
\usepackage{amsthm}
\usepackage{diagbox}
\usepackage{wrapfig}

\input pictex.tex

\theoremstyle{definition}

\input epsf

\begin{document}

\author{Nicol\'e Geyssel  \quad\quad  Mar\'ia Jos\'e Moreno \quad\quad Andr\'es Navas}

\title{On the geometry and topology of Da Vinci domes}
\maketitle

\noindent{\bf Abstract.} We study the famous Leonardo Da Vinci's domes, as well as the variations invented by Rinus Roelofs, from 
a mathematical viewpoint. In particular, we consider the problem of closing the dome in order to produce a spherical structure. We 
explain why this  problem is related to subtle geometric and topological considerations. This is in contrast with the 1-dimensional 
analog structure, namely Da Vinci's bridge, that can be easily closed up to make a circular shape. 

\vspace{0.2cm}

\section{Introduction}

Without any doubt, Leonardo Da Vinci is one of the great geniuses of history. His vast work is surprising not only for its depth, 
originality and beauty, but also for its heterogeneity. Among his writings, it is even possible to find several notes on mathematics. 
He got most of his mathematical education by his training under Luca Pacioli, for whom he illustrated the famous book {\em The 
Divine Proportion}. In this regard, the use of golden proportion in Da Vinci's art is frequently featured, yet usually exaggerated.

Da Vinci also made experiments --though not very substantive progress-- in other directions, for example, on the problems of squaring the circle and doubling 
the cube. For more information on this, see this article \cite{duv} of Sylvie Duvernoy. In particular, therein there is a nice discussion on a genuine and 
little known contribution of Da Vinci to the problem of determining the center of gravity of a pyramid.\footnote{In the literature, it is also mentioned  
that Da Vinci would have provided a very original proof of the Pythagorean theorem. However, this is highly unlikely, and the argument seems to have 
arisen much later. See this article \cite{pit} of Franz Lemmermeyer for a very complete discussion on this.}

In this article we will focus on a series of graphical sketches he left in his mythical work {\em Codex Atlanticus} from a mathematical 
(more precisely, geometric and topological) viewpoint.

\vspace{0.25cm}

\section{The Da Vinci bridge} 

Of the many structures Leonardo da Vinci designed, perhaps the most ingenious one is his design of an 
easy-to-assemble, self-supporting bridge. In ``the Codex'', he depicted this as follows:

\begin{center}
	\includegraphics[scale=0.34]{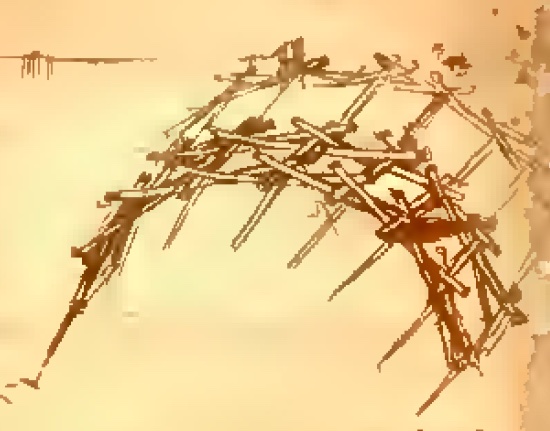}\\
	\scriptsize{Image from folio 69 recto of Leonardo Da Vinci's Codex Atlanticus}	
\end{center}

For a full description of the elementary geometry involved in this design, we refer to \cite{ham}. Over the last years, Da Vinci's 
idea has been systematically implemented as pedagogical activities to introduce students to 2D - 3D geometry in a more attractive 
way. (see for instance \cite{mariela}). It has also found a central place in mathematical exhibitions (museums, festivals, etc) 
throughout the world, and it has been beautifully displayed as decoration in many places, for example in Freiburg (Germany). 

\begin{center}
	\includegraphics[scale=0.16]{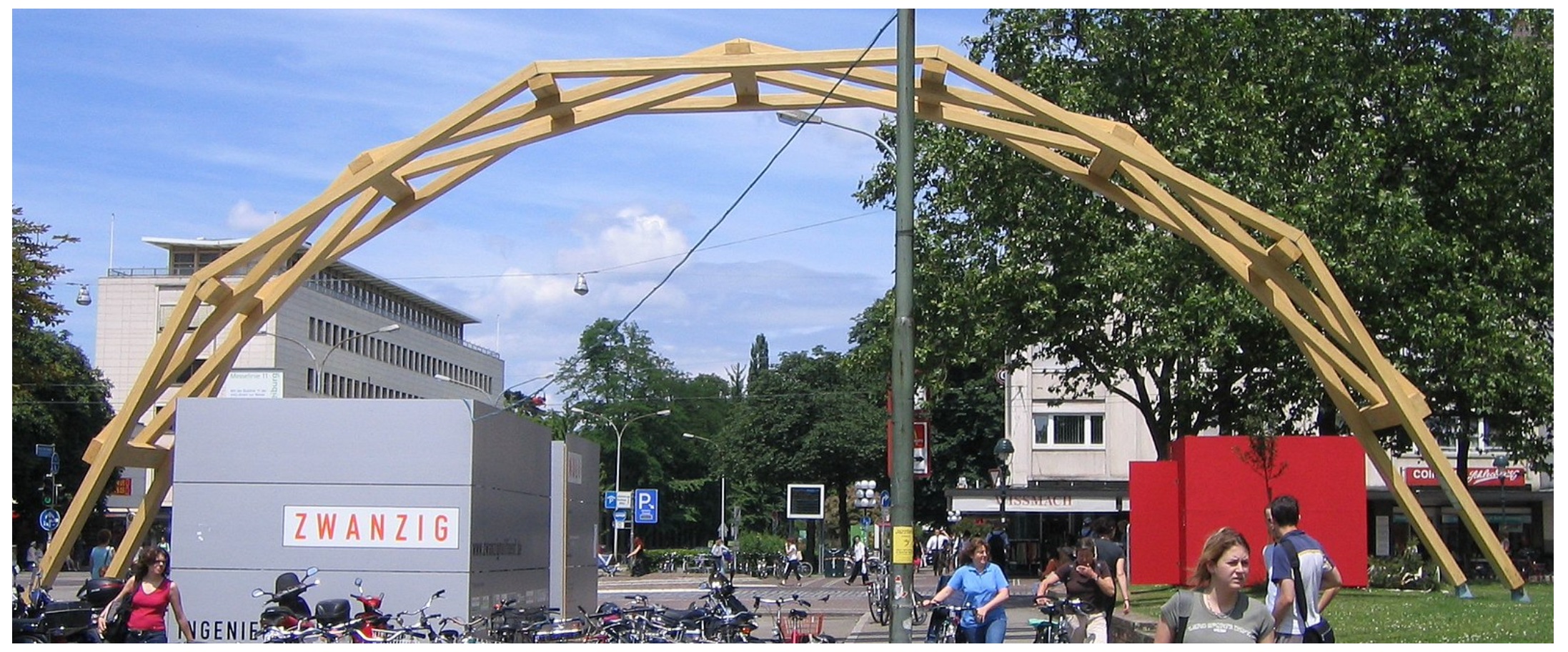}\\
	\scriptsize{A Leonardo Da Vinci's bridge in Freiburg; 
	image taken from Wikipedia commons}	
\end{center}

It is not very surprising that, by carefully pursuing this construction, in the end one obtains a ``circle'', as shown below.

\begin{center}
	\includegraphics[scale=0.26]{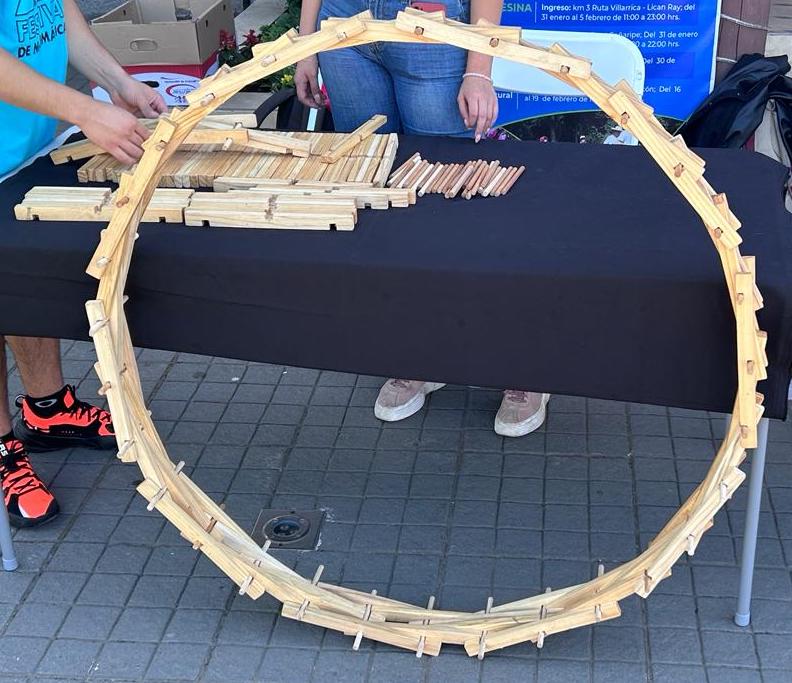}\\
	\scriptsize{A ``circle'' of Leonardo Da Vinci}	
\end{center}

\vspace{0.25cm}

\section{The Da Vinci domes}
\label{domes}

Here is another wonderful design from ``the Codex''. The structures arising from it are also self-supporting, 
but produce a kind of 2-dimensional surface instead of a 1-dimensional object. 

\begin{center}
	\includegraphics[scale=0.15]{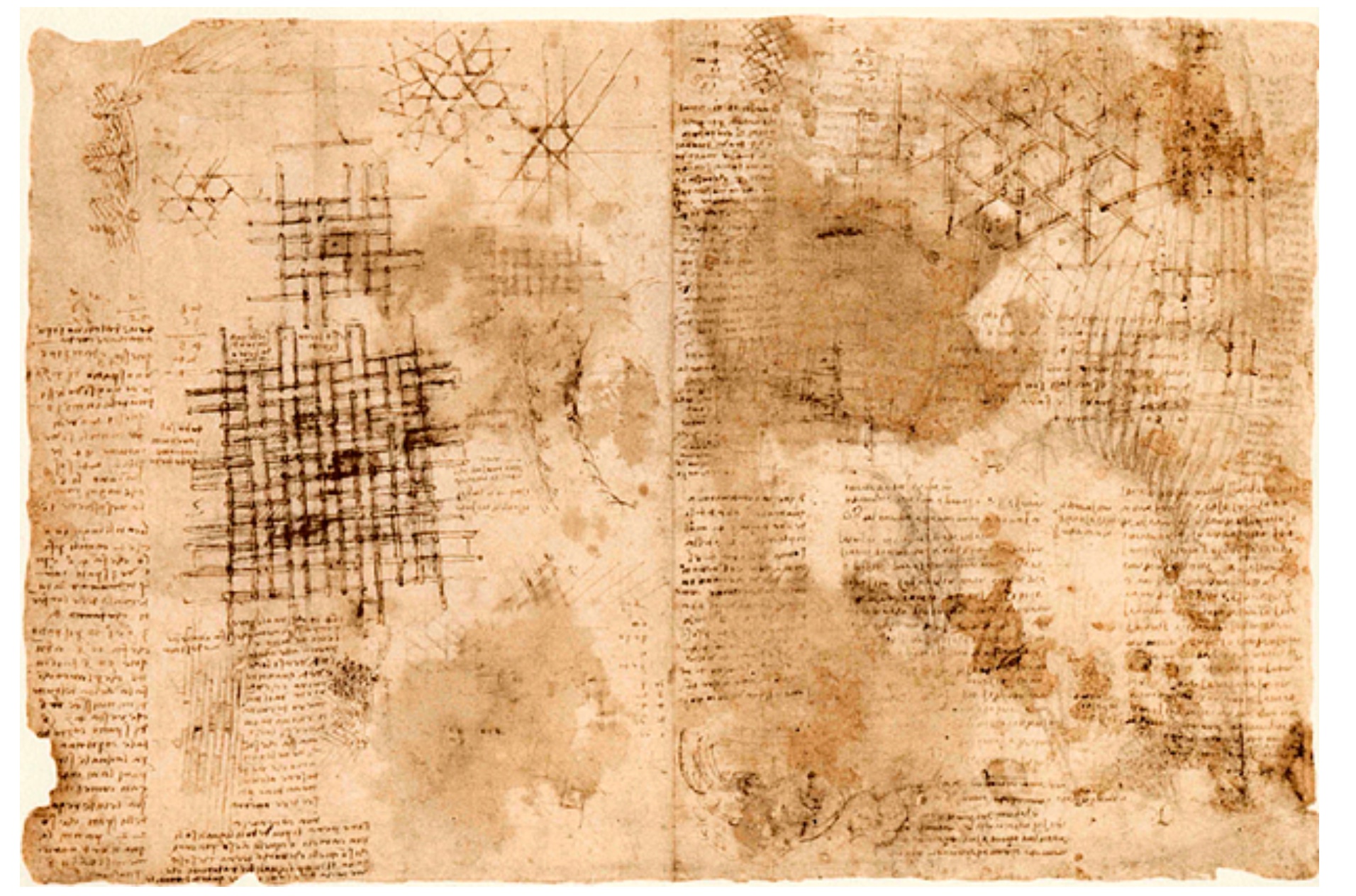}\\
	\scriptsize{Folio 899 verso of Leonardo Da Vinci's Codex Atlanticus}	
\end{center}

These structures are all based on the same principle: each rod is perfectly straight but has four notches. 
Quite crucially, the notches that are close to the ends are deeper than the other two. 

\begin{center}
	\includegraphics[scale=0.26]{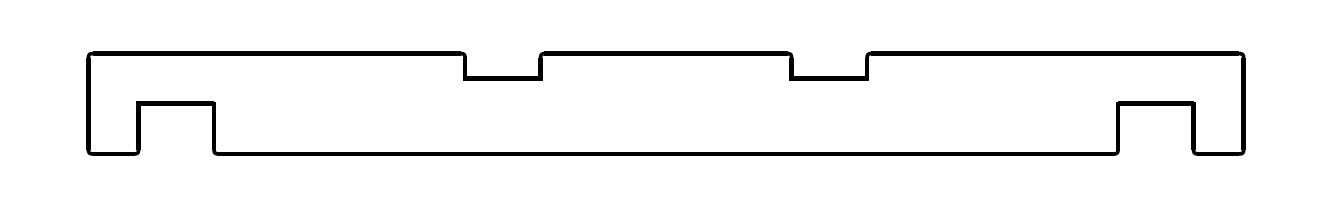}\\
	\scriptsize{A rod of a Da Vinci dome}	
\end{center}

When assembling pieces, each notch has to match with another one. 
Moreover, those in the ends have to be put towards the ground surface. Because 
of the difference of depth, this very simple rule produces elevation for the structure. 

\begin{center}
	\includegraphics[scale=0.16]{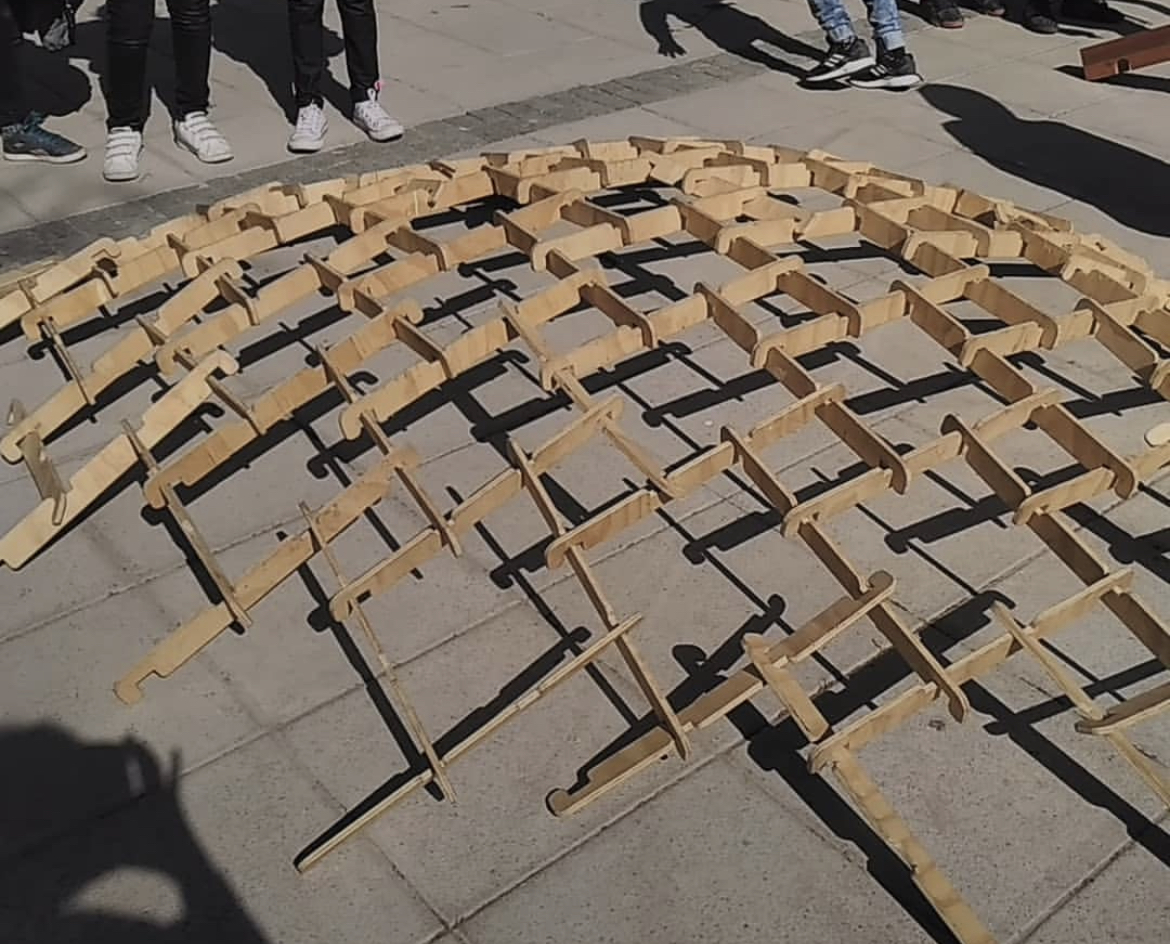}\\
	\scriptsize{A dome assembled for the Festival de Matem\'aticas in Valpara\'iso (Chile)}	
\end{center}

These structures have been rediscovered and pursued by Rinus Roelofs. We refer to \cite{rinus} for an account of this side of his 
amazing work. We will come back to this later. 

There  are many variations on the designs. In particular, the nice list of 11 periodic 
planar patterns below was produced for the LeonarDome educational  
kit  of the Museum of Mathematics of Catalonia; see \cite{catalunha}. 

\begin{center}
    \includegraphics[scale=0.57]{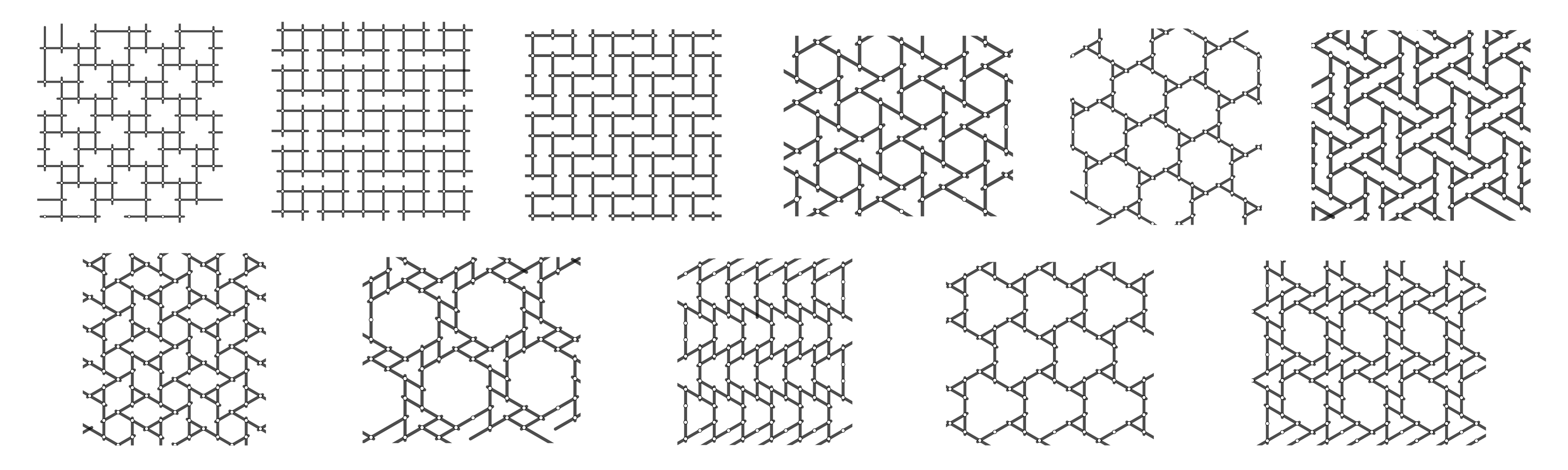}\\
	\scriptsize{Eleven different patterns to build domes}	
\end{center}

One can even create more designs. For example, this one was conceived by Ignacio (the son of the third-named author), 
yet we later discovered that it appears in (the top of) page 4 of \cite{forzado}.
\begin{center}
	\includegraphics[scale=0.06]{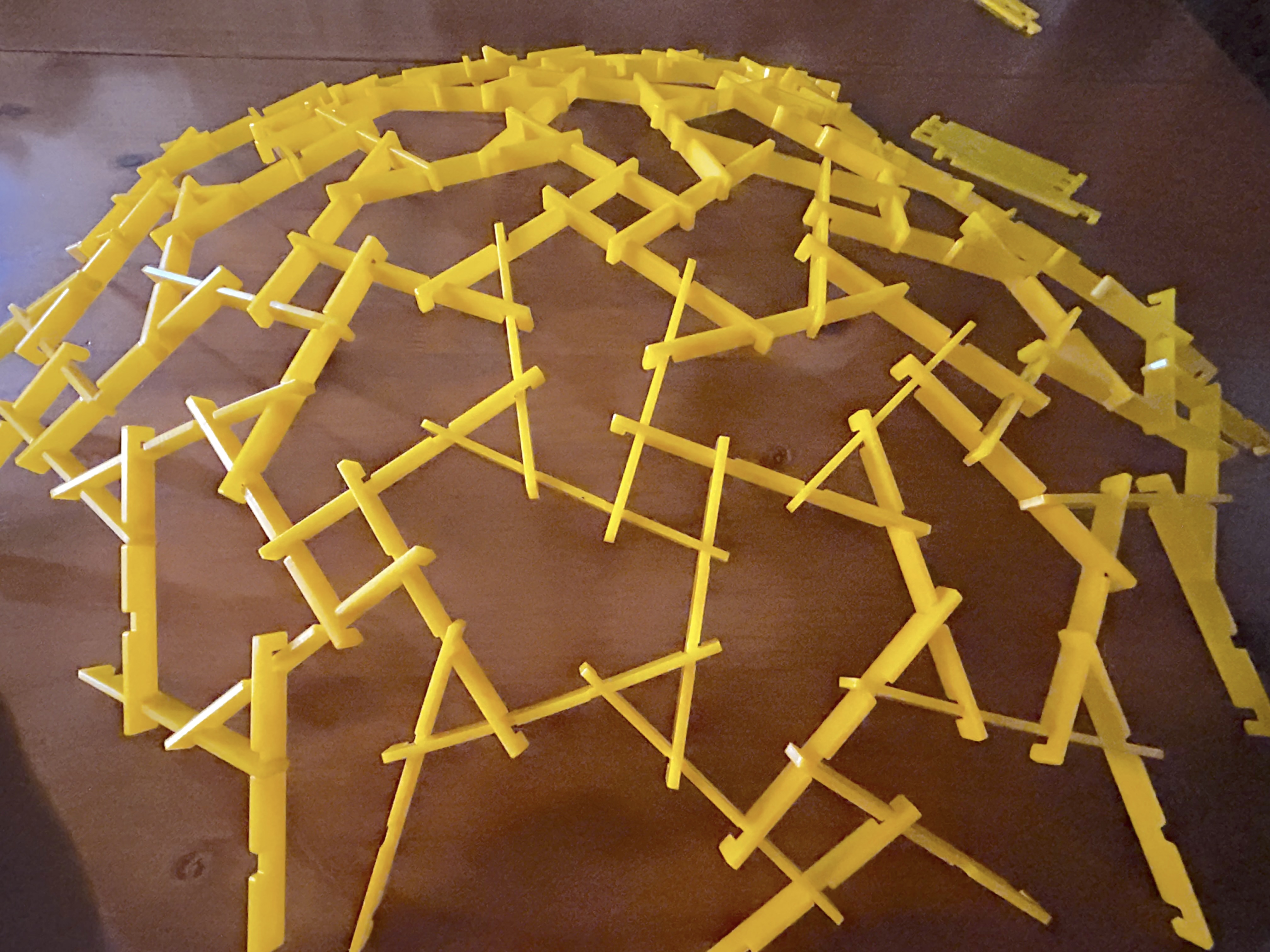}\\
	\scriptsize{A ``new'' pattern}	
\end{center}
We should stress that in the illustrations above, 
all edges have equal length. However, our discussion still applies to non-equal-length rods, which give raise to more patterns.

It seems interesting to explore the whole landscape of periodic planar tilings that one can implement this way.  
A challenging problem is to determine which of the 17 crystallographic groups can be realized. For instance, 
it seems that those that involve both nontrivial rotations and reflexions are hard to produce (some of these 
are perhaps impossible). One can readily check that the patterns above realize the groups pmg, pgg, p31m, 
p4g, p4 and p6 (see \cite{doris}  for a concise discussion on crystallographic groups, including notation).


\section{Da Vinci patterns and sphericity}

Da Vinci's bridge can be closed to a circular structure. Quite naturally, one may ask whether Da Vinci's domes can also be closed to spherical 
surfaces. It is worth stressing that we are only referring to the combinatorics of the patterns, with no claim about the physical stability of the 
structure. In particular, we do not suggest any tensegrity phenomenon here. 

Let us point out that there exist spherical tilings that respect combinatorial rules as those of Da Vinci's domes (interior notches 
of each rod point in one direction and those close to the ends in the opposite one). For instance, one can build them starting 
from ico-dodecaedral like configurations. The beautiful design below was conceived by Roelofs. Note, however, that the sphericity 
of the structure is somehow ``reverse'' to that of Da Vinci domes, and comes from the fact that the rods are systematically curved 
in a well-chosen direction (actually, there is no difference of depth between the notches).

\begin{center}
	\includegraphics[scale=0.32]{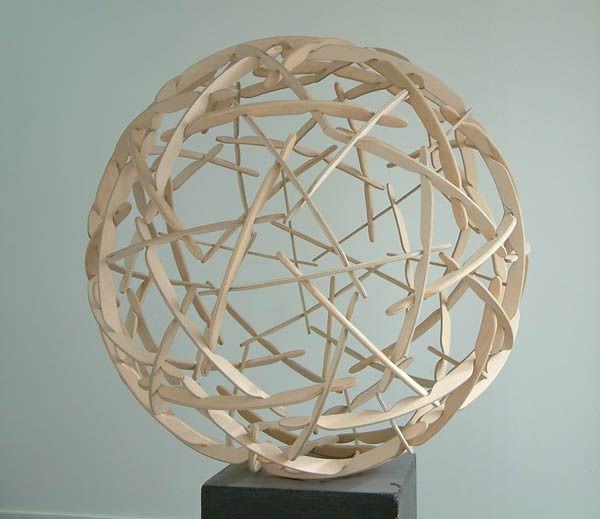}\\
	\scriptsize{Image taken from Roelofs' website: http://www.rinusroelofs.nl/structure/designs-lg/dm-sphere-10.html}	
\end{center}

Using  only straight rods still allows closing a spherical structure in certain cases. This is illustrated by the design below, 
due to Arie Brederode. Note that this requires a large number of rods. 
\begin{center}
	\includegraphics[scale=0.213]{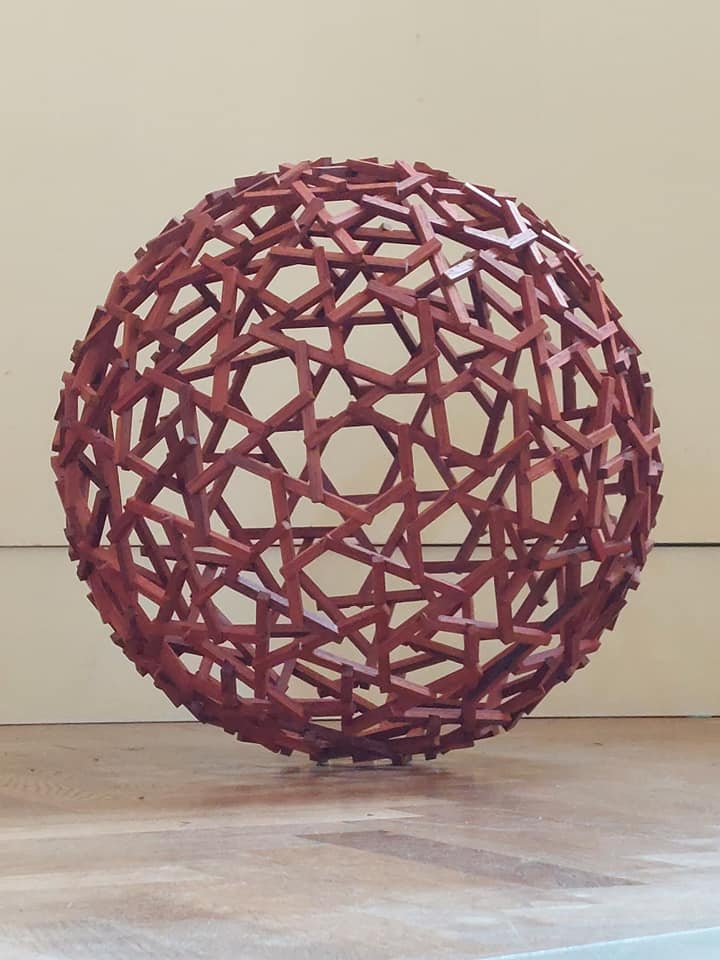} \qquad
        \includegraphics[scale=0.1]{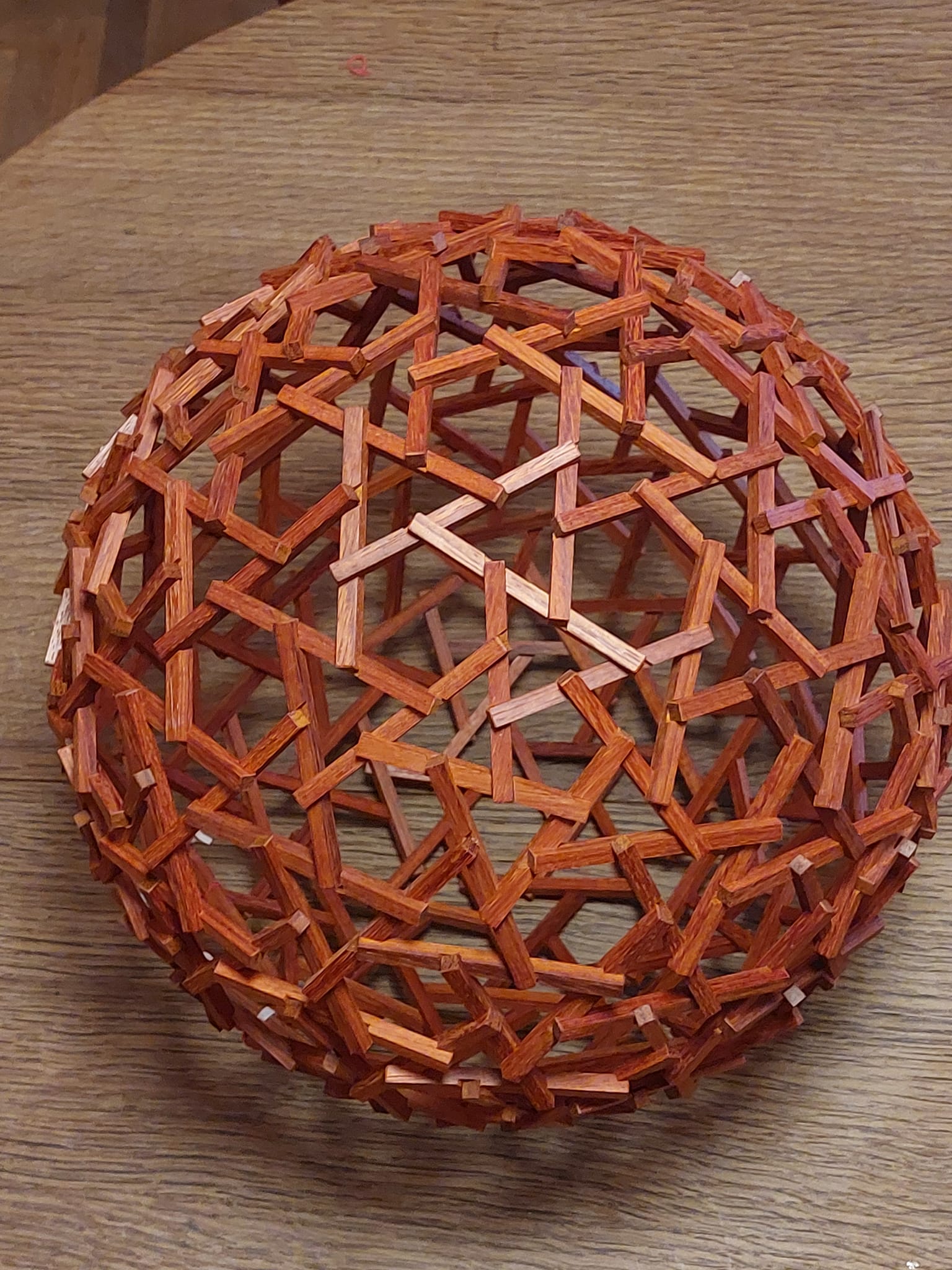}\\
	\scriptsize{Images uploaded by Arie Brederode to the Facebook page of Bridges: Arts \& Mathematics, reproduced with permission}	
\end{center}

\subsection{A warning on collinearity of vertices along rods}

From a theoretical viewpoint, modeling Da Vinci's domes made of straight rods as polyhedra in which vertices along the same rod are collinear 
is not appropriate (compare \cite{forzado}). Actually, such an ideal object can never close to a spherical surface. This is a consequence of the 
next elementary result, the hypothesis of which apply to Da Vinci domes if they were considered as polyhedra because of  the combinatorial 
constraint in their construction (more precisely, by the fact that no boundary notch of a rod matches with another boundary notch).

\vspace{0.38cm}

\noindent{\bf The Polyhedron Theorem.} {\em There is no spherical like polyhedron 
such that  two of the edges meeting at each vertex are contained in the same straight line.}

\vspace{0.38cm}

We provide two proofs of this result. The first one is purely geometric and somewhat connected to the idea of curvature. It 
reproduces a standard argument that shows that every compact surface in 3-dimensional space has points with positive 
Gaussian curvature. The second one, based on the classical Descartes angle defect theorem, has a topological flavor. Although 
it is more elaborate than the geometric one, the topological considerations involved will drive our further discussion. Quite 
naturally, both arguments proceed by contradiction.


\vspace{0.3cm}

\noindent{\bf First Proof (Da Vinci meets Gauss).}  
Assume by contradiction that a polyhedron as in the statement exists, and consider a very large sphere that completely contains it in the interior. 
By continuously  decreasing the radius of the sphere, we detect a moment where it becomes ``tangent'': the polyhedron touches the sphere, 
but doesn't cross it. Let $P$ be a point in common between the polyhedron and this new sphere. Such a point has to be an endpoint of any 
line that contains vertices of the polyhedron, otherwise this line would cross the sphere, and there would be vertices of the polyhedron 
outside the sphere. In particular, $P$ must be a vertex. However, by hypothesis, 
at each vertex of the polyhedron, there are at least two edges that are collinear, which implies that $P$ 
is not an endpoint of the corresponding line. This contradiction concludes the proof.

\vspace{0.3cm}

\noindent{\bf Second Proof (Da Vinci meets Descartes).}  Recall that the {\em angular defect} of a vertex polyhedron is defined as 
the difference between 360º and the sum of the angles at the vertex. An easy way to visualize this consists in mapping the faces 
converging to the vertex into a plane by opening one of the edges. The angle that remains in the plane uncovered by those of the 
faces is the angular defect.

\begin{center}
	\includegraphics[scale=0.13]{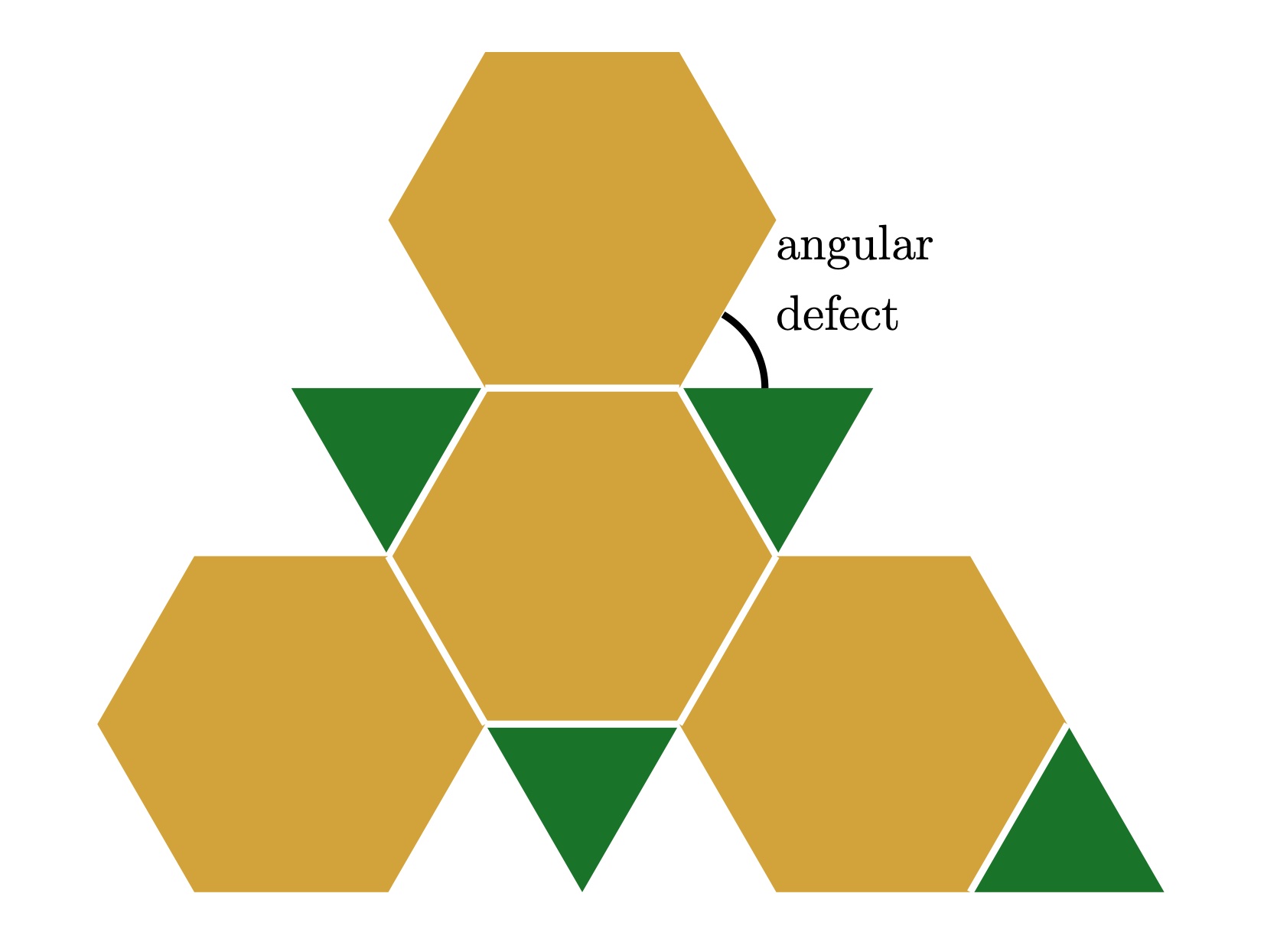}
\end{center}

A famous theorem due to Descartes establishes that, for a 
spherical like polyhedron, the total sum of the angular defects along the vertices equals 720º. (See \cite{descartes} for a proof 
of this theorem, which is an almost direct consequence of Euler's formula (\ref{euler}), yet it historically precedes it). 

Assume now that a polyhedron as in the statement of the theorem exists. If a vertex of such an object has degree 2, then the two angles 
converging to it are equal to 180º, hence the angular defect vanishes. If a vertex has degree 3, then it is shared by three angles: a pair of 
supplementary ones and another one of 180º. The sum of these angles would hence equal 360º and, therefore, the angular defect would 
also vanish. In case of higher degree, one can still see that the angular defect is nonpositive. Therefore, the total sum of the angular defects 
along vertices is nonpositive. However, this contradicts Descartes' theorem.

\begin{center}
	\includegraphics[scale=0.25]{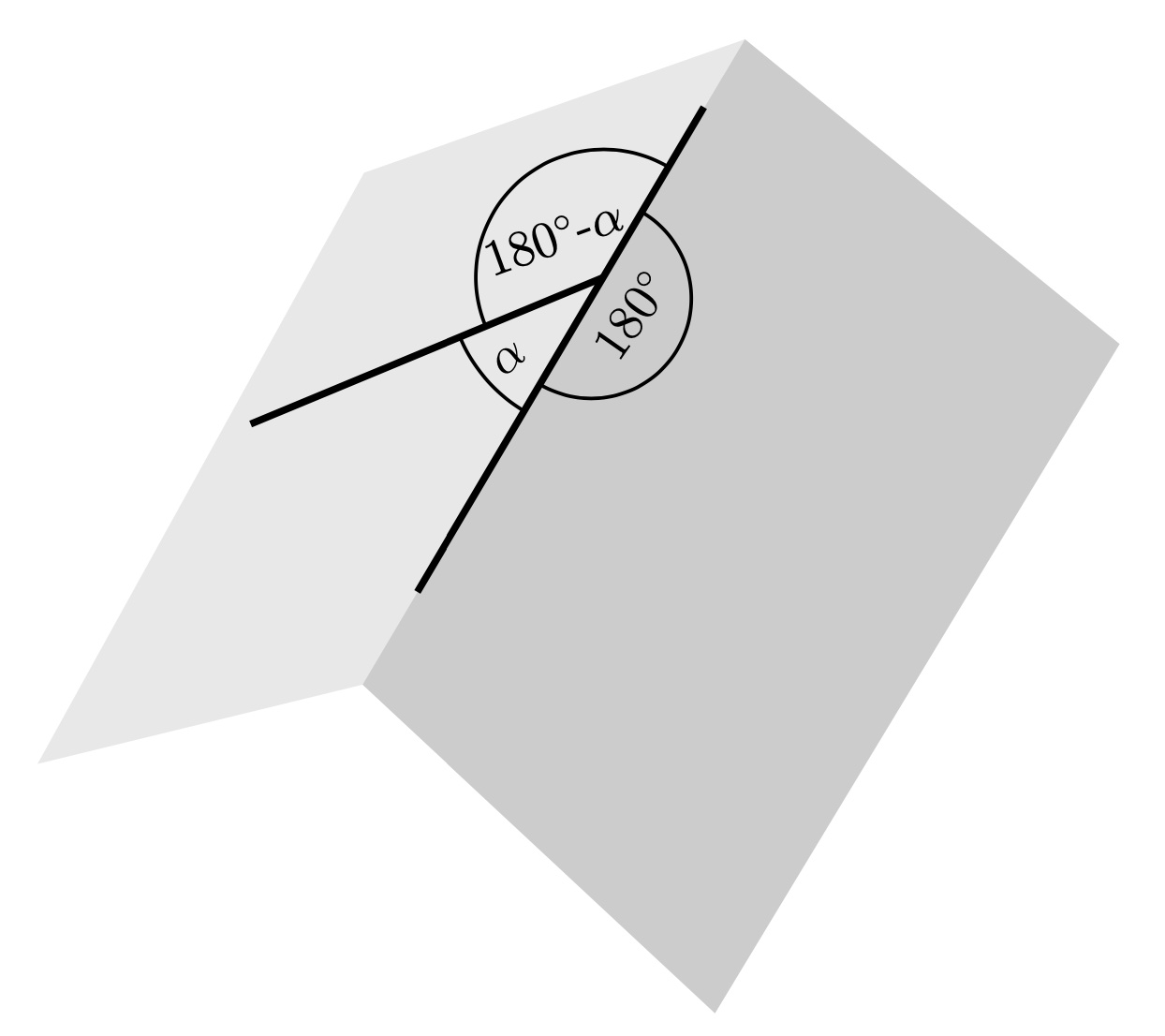}\\
	\scriptsize{Vanishing of the angular defect at each degree-3 vertex}	
\end{center}

One can argue that the previous theorem is badly suited because allowing only collinear vertices along the same rod doesn't imply that 
the whole object can be considered as a polyhedron, since the ``faces'' could be nonplanar. However, if we artificially add edges to the 
object thus produced so that it becomes a genuine polyhedron (one can easily do this by triangulating the faces), the theorem above still 
applies (this is the reason why we included the possibility of arbitrary degree for vertices in the statement). 

The theorem above emphasizes once more the relevance of the difference of depth between interior and exterior notches of the rods. This 
allows producing angular defect (after triangulation, so that the object in view is a genuine polyhedron) which, as we have seen, is necessary  
to close a spherical structure. Since the angular defect thus produced is very small, a large number of vertices is needed to complete the 
total angle 720º predicted by Descartes theorem. This is the reason why a large amount of edges is needed. 

For the case of curved rods, part of the gaining of angular defect comes from the curvature of the rod. Of course, in order to make this 
claim formal, one previously needs to forget about the rods keeping in mind only the vertices, and connect them by straight lines. Doing this, 
one will immediately notice that vertices corresponding to the same rod are far from being collinear. This is the reason why so few rods are 
used in Roelof's sculpture. However, as we will see in the next section, the careful choice of the combinatorial structure is also crucial. 


\section{Da Vinci meets Euler: exploring topology}

Although certain Da Vinci like structures can be closed to make a spherical one, this is impossible with the eleven patters previously depicted, 
as well as with all patterns coming from periodic planar tilings. The reason relies on Euler type theorems, as we next explain.

Suppose that we have a spherical like structure, and project it from inside to a sphere surrounding it. We thus obtain a spherical graph for 
which we have a number of vertices, of curves connecting them that we still call ``edges'', and regions limited by these edges that we still call 
``faces''. Denote by $V$, $E$ and $F$ the number of these vertices, edges and faces, respectively. Euler's theorem establishes the equality 
\begin{equation}\label{euler}
 V - E + F = 2
 \end{equation}
(see \cite{erdos} for an elementary proof; see also \cite{yo} for a discussion in relation to designs over spheres).  
Now, each pattern depicted in \S \ref{domes} exhibits translational symmetries in two independent directions. 
By well-known facts of topology, this implies that, if it closes a surface, then  
\begin{equation}\label{torus}
V-E+F = 0,
\end{equation}
which is in contradiction with (\ref{euler}). 

Here is a concrete example. The pattern below involves triangles and dodecagons (in order to apply Euler type formulae, we need to 
count a vertex at each crossing of rods). Let $T$ and $D$ denote the number of the former and the latter, respectively. 
Since each vertex involves three edges, we have
$$V = \frac{3 \, T + 12 \,D}{3}.$$
Each edge is shared by two faces, hence, 
$$E = \frac{3 \, T + 12 \,D}{2}.$$
Finally, the number of faces is the sum of the number of triangles and dodecagons: 
$$F = T + D.$$ 
Therefore, 
$$V - E + F = \frac{3 \, T + 12 \,D}{3} - \frac{3 \, T + 12 \,D}{2} + (T + D) = \frac{T}{2} - D.$$
However, one can directly check by looking at the pattern that:

\vspace{0.1cm}

\noindent -- each dodecagon is neighbored by six triangles;

\vspace{0.1cm}

\noindent -- each triangle is surrounded by three dodecagons.

\vspace{0.1cm}

\noindent This immediately implies that $T = 2D$, hence 
$$V-E+F = \frac{T}{2} - D = 0.$$

\begin{center}
	\includegraphics[scale=0.65]{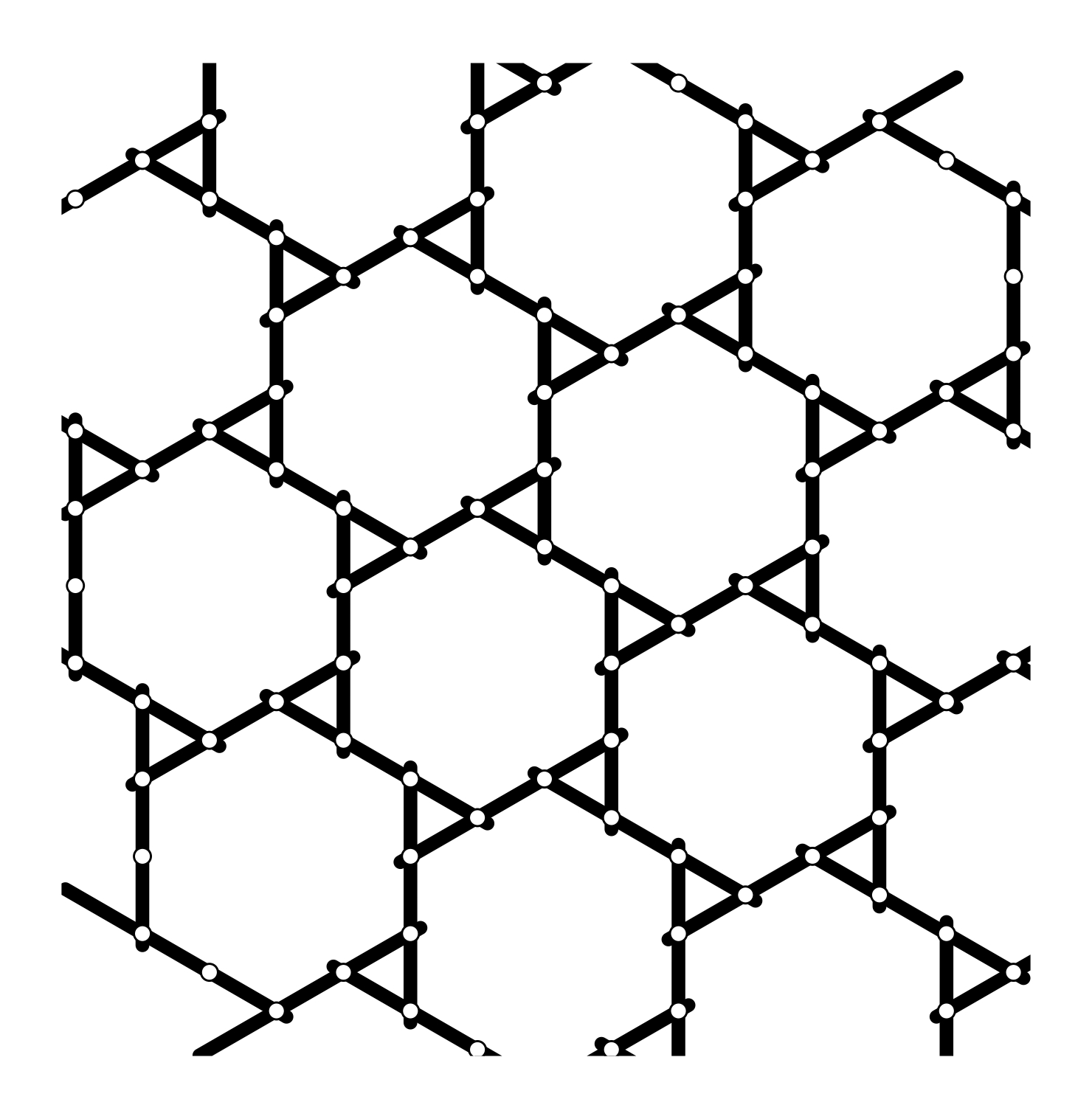}\\
	\scriptsize{Proving the equality $V-E+F=0$ for a concrete example \\ 
	(the reader is invited to give similar arguments for each of the eleven patterns previously depicted)}	
\end{center}

\noindent{\bf A technical remark.} 
There is a way of showing (\ref{torus}) for all planar patterns above that only uses the (more elementary) Euler's formula (\ref{euler}) 
of spherical tilings. This comes from the fact that these patterns yield tilings on the torus, for which (\ref{torus}) holds. Roughly, the proof 
proceeds as follows. Using the two independent symmetry directions, close the pattern so as to build a torus surface, and 
let $V,E,F$ be the number of vertices, edges and faces, respectively. 
Now, coming back to the plane, copy the planar patterns 
$n$ times horizontally and vertically, and denote $V_n, E_n, F_n$ the number of vertices, edges and faces that do arise. Non-closeness only 
happens along the boundary. Since this boundary ``grows linearly'', the differences $|V_n - n^2 V|$, $|E_n - n^2 E|$ and $|F_n - n^2 F|$ 
are bounded by a number that is linear in $n$. All of this implies that 
$$V_n - E_n + F_n = n^2 \, (V-E+F) + O (n).$$
Now, this $n^2$-times copied pattern induces a (non Da Vincian !) 
tiling of the sphere just by adding a ``face towards infinity''. By Euler's formula (\ref{euler}) for 
general spherical tilings, we have 
$$V_n - E_n + (F_n+1) = 2.$$
Therefore, 
$$n^2 \, (V-E+F) = O (n).$$
Dividing by $n^2$ and letting $n \to \infty$, this shows the announced equality $V-E+F=0$.


\section{Conclusions}

The previous discussion shows that, in order to pursue with new Da Vinci like spherical designs, it is worthwhile to take into account not only 
geometric, but also topological constraints. Roelofs design is an example of this: it is based on (actually, combinatorially equivalent to) a 
truncated icosahedron. It seems possible to use other classical designs to produce Da Vinci like versions; a challenging problem would be 
to realize models inspired from different Archimedian solids. Note that Brederode's design also arises from a well-known family of polyhedra, 
namely Goldberg's polyhedra.

Needless to say, besides geometric and topological issues, physical considerations of the material employed  for the concrete realizations 
also needs to be taken into account. For example, the width of the rods creates problems at the junctures to produce polygonal 
configurations of a very large number of sides. 

We stress that our considerations do not close the study of designs that yield to other kind of surfaces. The first proof of the Polyhedron Theorem  
applies in more generality to closed like polyhedra, as for instance toroidal ones. However, considering notches of different depth could lead 
to interesting structures of this kind. Note that there is no combinatorial obstruction for planar configurations to yield toroidal structures, since 
a torus is obtained by a well-known process of ``folding a piece of paper in two directions''. As an example we produced the figure below, which 
shows that the first of the 11 plane tilings can be implemented on a toroidal like surface. Producing this way a more symmetrical and attractive  
object would certainly require the work of an artist.

\vspace{0.2cm}

\begin{center}
	\includegraphics[scale=0.26]{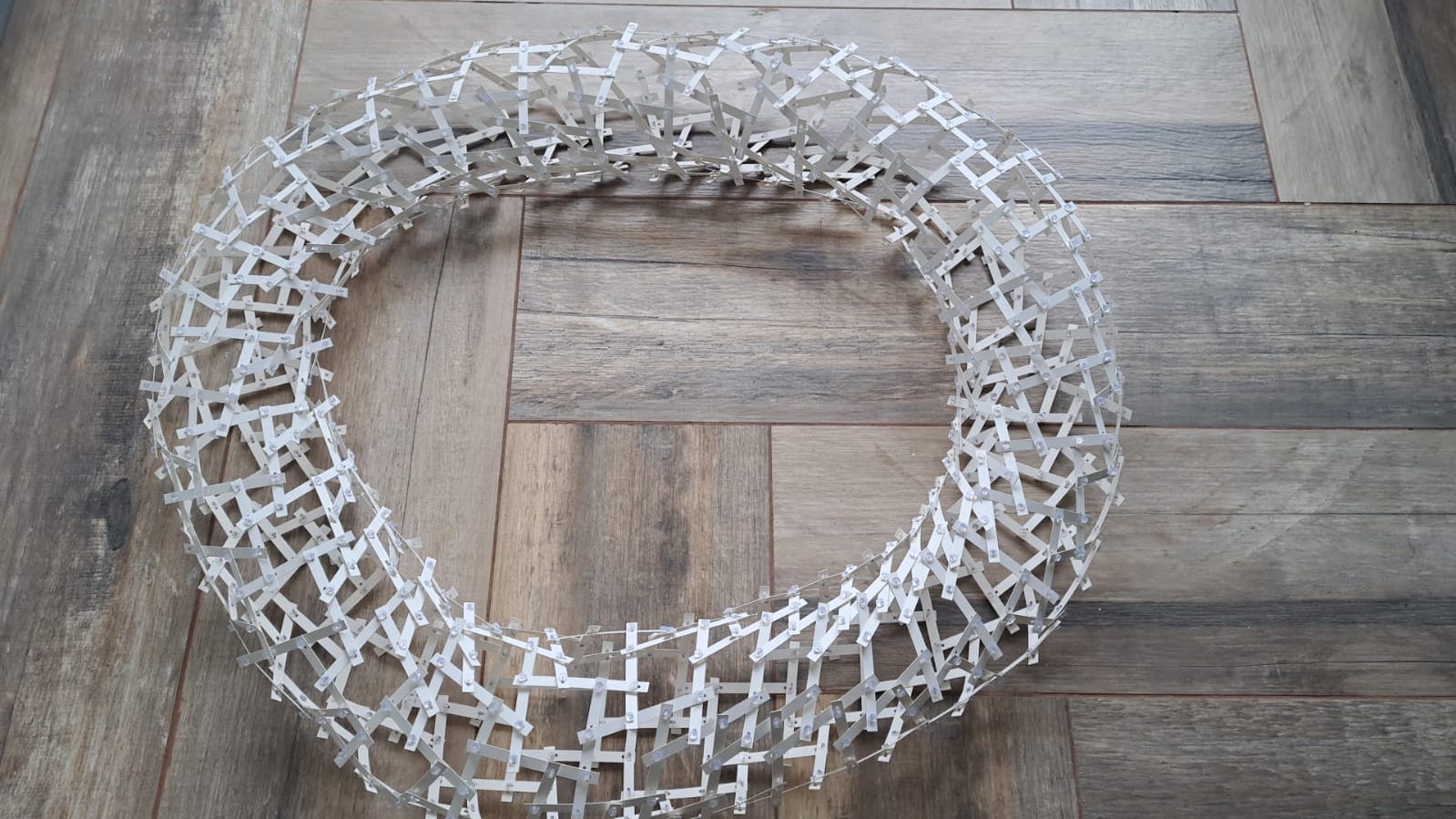}\\
	\scriptsize{One of the eleven plane tilings implemented on a torus}	
\end{center}

\vspace{0.2cm}

\noindent{\bf Acknowledgments:} We would like to thank Mario Ponce for his interest on this text and his many remarks to it.


\begin{small}

\vspace{0.1cm}

\noindent Nicol\'e Geyssel, UTFSM, El Pinar 36, San Joaqu\'in, Santiago, Chile (nicole.geyssel@usm.cl)\\

\vspace{0.1cm}

\noindent Mar\'ia Jos\'e Moreno, PUC, Vicuña Mackena 4860, Macul, Santiago, Chile (mmors@uc.cl)\\

\vspace{0.1cm}

\noindent Andr\'es Navas, USACH, Alameda 3363, Est. Central, Santiago, Chile (andres.navas@usach.cl)\\

\end{small}


\begin{thebibliography}{Dillo 83}

\bibitem{erdos} 
{\sc Martin Aigner \& G\"unter M. Ziegler.} 
{\em Proofs from The Book.} Springer Verlag (1998).

\bibitem{catalunha} 
{\sc Enric Bras\'o.}
Les c\'upules de Leonardo da Vinci. 
{\em Noubiaix} {\bf 42} (2018), 116-126.

\bibitem{mariela} 
{\sc Pilar Carrasco, Mariela Carvacho \& Francisco S\'anchez.}
C\'upulas de Da Vinci: la geometr\'ia colaborando con otras disciplinas. 
En {\em Educaci\'on Matem\'atica Insterdisciplinar en el Aula}, by 
Jaime Huinacahue and Daniela Soto eds., Ediciones UCM (2023).

\bibitem{duv} 
{\sc Sylvie Duvernoy.} 
Leonardo and theoretical mathematics. Leonardo da Vinci: architecture and mathematics.  
{\em Nexus Netw. J.} {\bf 10}, No. {\bf 1}, (2008), 39-49.

\bibitem{ham}
{\sc Hans Humenberger.}
Mathematische Aktivitäten rund um die Leonardo-Brücke. 
In: {\em Schriftenreihe zur Didaktik der Mathematik der Österreichischen Mathematischen Gesellschaft} 
{\bf 45} (2021), 67-86. 

\bibitem{pit} {\sc Franz Lemmermeyer.} 
Leonardo da Vinci's proof of the Pythagorean Theorem. 
{\em The College Mathematics Journal} {\bf 47} (2016), 361-362. 

\bibitem{yo} {\sc Andr\'es Navas.} 
Une \begin{tiny}$^{_\ll}$\end{tiny}erreur\hspace{0.04cm}\begin{tiny}$^{_\gg}$\end{tiny}  
 g\'eom\'etrique dans la Ligue des Champions. 
{\em Images des Mathématiques}, CNRS (2019). 

\bibitem{rinus} {\sc Rinus Roelofs.} 
Two- and three-dimensional constructions based on Leonardo grids. 
{\em Nexus Netw. J.} {\bf 10}, No. {\bf 1} (2008), 17-26.

\bibitem{doris} {\sc Doris Schattschneider.} 
The plane symmetry groups: their recognition and notation. 
{\em The American Mathematical Monthly} {\bf 85} (1978), 439-450.

\bibitem{descartes} {\sc Paul Scott.} 
Angle defect and Descartes' theorem.
{\em Australian Mathematics Teacher} {\bf 62}, nº {\bf 1} (2006), 2-4.

\bibitem{forzado} {\sc Peng Song, Chi-Wing Fu, Prashant Goswami, Jianmin Zheng, Niloy J. Mitra \& Daniel Cohen-Or.} 
Reciprocal Frame Structures Made Easy. 
{\em ACM Transactions on Graphics} {\bf 32}, Issue {\bf 4}, Article No. {\bf 94} (2013), 1-13.

\end{thebibliography}
\end{document}